\documentclass[10pt]{article}
\usepackage{amssymb}
\usepackage{amsmath,amsfonts, verbatim,ifthen}
\usepackage{amsthm}
\usepackage{color}
\usepackage{multirow} 
\usepackage[all]{xy}
\usepackage{graphicx}
\DeclareGraphicsExtensions{.pdf,.png,.jpg}

\renewenvironment{proof}[1][Proof]{\textbf{#1} }{\ \rule{0.5em}{0.5em}}

\newcommand{\ignore}[1]{}

\numberwithin{equation}{section}

\begin{document}

\title{On the Additive Volume of Sets of Integers}
\author{Gregory A. Freiman \\ School of Mathematical Sciences,\\
Tel Aviv University, Tel Aviv 69978, Israel. \\
e-mail:
\textsf{grisha@post.tau.ac.il}.}

\maketitle

\thispagestyle{empty}

\lineskip=2pt\baselineskip=5pt\lineskiplimit=0pt

\noindent{\textsc{Abstract}}:  \begin{quote} The study of the additive volume of sets can be reduced to the case of one-dimensional sets. The exact values of the volume of extremal sets are given as a conjecture.
\end{quote}

\bigskip
\noindent {\emph{AMS} classification
numbers: 11}

\bigskip


\newpage
\lineskip=1.8pt\baselineskip=18pt\lineskiplimit=0pt \count0=1

The exposition in this paper is neither complete, nor perfected. Nevertheless, I decided to publish it because of the general interest and the importance of the theme and the new ways I propose to treat it.

\section{Introduction}

Our main object of study is a finite set $A\subset \mathbb Z$. We say that  $A$ is in normal form when 
\[
A=\{a_0=0<a_1<\cdots<a_{k-1}\} 
\]
and ${\rm g.c.d.}(a_1,\dots,a_{k-1})=1$.

We call {\it doubling\/} of $A$ the set
\[
A+A:=2A= \{x:\, x=a+b,\ a,b\in A\}.
\]

Our exposition will use two basic numbers characterising the set $A$:

the cardinality of $A$
\[
k:=|A|;
\]

the cardinality of the doubling of $A$
\[
T:=|2A|.
\]

The main notions we will use are those of Freiman homomorphism and isomorphism,
dimension, parallelepiped, $n$-dimensional arithmetic progression, and the additive volume of the set $A$.

The main previous result relevant to us is Freiman's Theorem, see Nathanson [13].

The purpose of this paper is to show how the study of sets of arbitrary dimension can be carried out most effectively by transition to the case of one-dimensional sets.

For the convenience of readers who are not familiar with the above notions we will give now the main definitions.

{\bf Definition 1} (Nathanson [13],  p. 233) Let $G$ and $H$  be Abelian groups and let $A\subseteq G$ and $B\subseteq H$. A map $\phi:A\to B$ is called a {\it Freiman homomorphism\/}
of order 2 if
\begin{equation}
\phi(a_1)+\phi(a_2)=\phi(a_1')+\phi(a_2')
\end{equation}
for all $a_1,a_2,a_1',a_2'\in A$ such that
\begin{equation}
a_1+a_2=a_1'+a_2'.
\end{equation}

In this case the induced map $\phi^{(2)}:2A\to 2B$, given by
\begin{equation}
\phi^{(2)}(a_1+a_2)=\phi(a_1)+\phi(a_2),
\end{equation}
is  well defined. 

If $\phi: A\to B$ is a one-to-one correspondence such that (3) holds if and only if (2) holds, then $\phi$ is a {\it Freiman isomorphism\/} of order 2 and the induced map $\phi^{(2)}:2A\to 2B$ is also one-to-one.

In this paper we will simply call these notions homomorphism and isomorphism.

{\bf Definition 2.} A {\it parallelepiped\/} is a set
\[
D=(a_1,\dots,a_n)+\{(x_1,\dots,x_n)\in\mathbb Z^n:\, 0\le x_i< h_i,\,\, i=1,\dots,n\},
\]
where $h_i\in\mathbb Z$, $i=1,\dots,n$.

I. Rusza introduced the useful notion of $n$-dimensional arithmetic progression.

{\bf Definition 3.} Let $a,q_1,\ldots,q_n$ be elements  of $\mathbb Z$, and let $l_1,\ldots,l_n$ be positive integers. The set
\[
Q=Q(a,q_1,\ldots,q_n;l_1,\ldots,l_n)=\{a+x_1q_1+\cdots+x_nq_n:\, 0\le x_i< l_i\,\,{\rm for}\,\, i=1,\ldots,n\}
\]
is called an {\it $n$-dimensional arithmetic progression\/} in the group $\mathbb Z$. We usually take $a=(a_1,\ldots,a_n)$ $=0$.

The length of $Q$ is defined to be $l(Q)=l_1 \cdots l_n$. Clearly, if we denote the cardinality of a set $X$ by $|X|$, then
\[
|Q|\le l(Q).
\]

It is clear that an $n$-dimensional arithmetic progression is the image of a suitable $n$-dimensional parallelepiped under a homomorphism. If the latter is an isomorphism, the progression is said to be {\it proper\/}.

We quote  Melvyn B. Nathanson (see [13], p. 231). 

``Freiman's theorem asserts that if $A$ is a finite set of integers such that the sumset $2A$ is small, then $A$ is a large subset of a multidimensional arithmetic progression''.

More precisely,  the following results holds.

{\bf Theorem 1.} (Freiman) {\it Let $A$ be a finite set of integers such that
\begin{equation}
|2A|\le c|A|.
\end{equation} Then there exists integers $a,q_1,\ldots,q_n;l_1,\ldots,l_n$ such that, if $Q$ is as in Definition $3$, then $A\subseteq Q$,
where 
\begin{equation}
|Q|\le c'|A|
\end{equation}
and $n$ and $c'$ depend only on $c$.} 

After  this theorem, many results were obtained, gradually improving the estimate of the  constant $c'=c'(c)$. For instance:
\vspace{2mm}

Y.  Bilu [1]:
\[ 
c'  \le (2d)^{\exp\exp\exp(9c\,\log(2c))}.
\]

I. Rusza [14]:
\[
c'\le \exp\{ c^{c^c}\}
\] 

Mei Chu-Chang [2]:
\[
 c'\le \exp\{Kc^2(\log c)^3\}, 
\]
\quad\,\ where $K$ is a constant.

T. Sanders [15]:
\[
c'\le \exp\{O(c^{7/4}\log^3c)\}.
\]

S. V. Konyagin [11]:
\[
c'\le \exp\{c\log c\}.
\]
 
T. Schoen [16 ]:
\[
c'\le \exp\{c^{1+D(\log c)^{-1/2}}\},
\]
where $D$ is a positive constant.

The notion of  isomorphism enables us now to introduce the notion of   dimension of a set $A$.

{\bf Definition 4.} We say that a set $A\subset\mathbb Z^n$  has {\it dimension\/} $d$ if there exists a set $B\subset\mathbb Z^d$ isomorphic  to $A$ such that $B$ is not contained in any hyperplane of $\mathbb Z^d$, and $d$ is the maximal number with this property. We denote it by ${\rm dim}\, A$.

We are now ready to introduce the notion of additive volume of a set $A$ with given $|2A|$.

{\bf Definition 5.} Let $A$ have dimension $d$. Among the sets $B$ as in the preceding definition, take a $B$ such that the number $V(B)$ of integer points in the convex hull of $B$ is minimal. We call this number $V(B)$   the {\it additive volume} of the set $A$, and denote it by $V(A)$, observing immediately that it is simultaneously the volume of any other set isomorphic to $A$.

Let us now formulate the main idea of the present paper. Two sets $A$ and $A_1$ with the same
characteristics $k$ and $T$, i.e., such that $k=|A|=|A_1|$ and $T=|2A|=|2A_1|$ are not at all necessarily isomorphic, and their volumes $V(A)$ and $V(A_1)$ do not have to be equal.  One can ask if, generally, is it possible, instead of a given family of isomorphic sets, to consider another family, with the same or close characteristic values $k$ and $T$, and which is simpler to study? One of our main results reads:
 
{\bf Theorem 2.} {\it Let  ${\rm dim}\, A=2$. Then there exist a set $A_1\subset \mathbb Z$ and a one-to-one homomorphism $\phi:A\to A_1$  such that 
\[
  {\rm dim}\, A_1=1,
\]
\[
T(A_1)<T(A),
\]
and
\[
V(A_1)>V(A).
\]
}

\section{From dimension 2 to dimension 1}

{\bf 2.1.}  When we compute $T=|2A|$ we have to take into account relations of the form
\begin{equation}
a_i+a_j=2a_s
\end{equation}
and
\begin{equation}
a_i+a_j=a_r+a_s
\end{equation}
among the elements of $2A$.

S. Konyagin and V. Lev [12] found a very useful representation of  relations (2.1) and (2.2) by means of $k$-dimensional integer vectors. Let us explain it on the following example, which makes clear the general case. Take $k=7$ and consider the relations
\begin{equation}
a_0+a_2=2a_1
\end{equation}
and
\begin{equation}
a_1+a_5=a_2+a_4.
\end{equation}
We associate to (2.3) and (2.4) the 7-dimensional vectors
\[
(1,-2,1,0,0,0,0)
\]
and
\[
(0,1,-1,0,-1,1,0),
\]
respectively. 

In the general case a system of relations (2.1) and (2.2) similarly leads to a collection of $k$-dimensional integer vectors. We can now look for a maximal linearly independent system of vectors in this collection, to which corresponds a maximal   linearly independent system of relations; the cardinality of this system of vectors (or relations) is denoted by $\lambda(A)$.

{\bf Theorem 3.} (Konyagin and Lev) \emph{For any   set $A\subseteq \mathbb Z^n$ of $k$ elements,  the dimension of $A$ is equal to
\[
{\rm dim}\,  A=k-1-\lambda(A).
\]
}

We can use the above connection and tools of linear algebra in order to compute the number $\lambda(A)$ and find explicitly a maximal system of linearly independent relations.

{\bf 2.2.}  Let $A$ be a set for which there is $D(A)$ such that $\bar a=0$ and $h_2=2$, i.e., $A$ lies in the union of the two lines $y=0$ and $y=1$, and $0\le x\le h_1-1$.

Which points in $D(A)$ are in $A$?

Let $A'=\{\bar z:=(x,y): \, y=0\}$ and $A''=\{\bar z=(x,y):\, y=1\}$.
The sets $A'$ and $A''$ can be arbitrarily moved along the axis $x$. So, let us assume that $(0,0)$ and $(0,1)$ lie in $A$. Suppose also that $(x_0,0)$ and  $(x_1,1)$ are in $D(A)$, and $x_0$ and $x_1$ are maximal: they are the right end-points of the segments that contain $A$. One can consider that $x_0\ge x_1$. Therefore, $(h_1-1,0)\in A$.
\begin{equation}
\overset{\displaystyle{(0,1)}}{\overset{\displaystyle\bullet}{\overset{\displaystyle\vert}{\underset{\displaystyle{(0,0)}}\bullet}}}-------------
\underset{\displaystyle{(h_1-1,0)}}{\bullet}
\end{equation}

Thus, summarizing, the points $(0,0)$, $(0,1)$, and $(h_1-1,0)$ lie in $A$, and the point $(b,1)$ with $b\le h_1-1$ is the extremal right point on the line $y=1$. It can even coincide with $(0,1)$,

Let us project $A''$ on the line $y=0$ parallel to the vector $\ell=(\ell_1,-1)$, where first $\ell_1> 2h_1-2$, and denote the projection map by $\phi$:
\begin{equation}
\underset{\displaystyle{(0,0)}}\bullet-------
\underset{\displaystyle{(h_1-1,0)}}\bullet-------
\underset{\displaystyle{(\ell_1,0)}}\bullet-------\underset{\displaystyle{(\ell_1+b,0)}}\bullet
\end{equation}
Then $A_{\rm proj}$ is isomorphic to $A$.

Now let us put
\[
\ell_1=2h_1-2.
\]

The matrix of the map $\phi$ is
\[
M=\left( \begin{array}{cc}
1 & 0  \\
2(h_1-1) & 0\end{array} \right),
\] 
so that
\[
(x,y)\left( \begin{array}{cc}
1 & 0  \\
2(h_1-1) & 0\end{array} \right) =(x+2(h_1-1)y,0);
\]
note that 
\[
\phi(x,1)=(x+2(h_1-1),0).
\]

We have $A_1:=A_{\rm proj}=A'\cup A''_{\rm proj}$,
where
\[
A'\subset [0,h_1-1],
\]
\[
A''_{\rm proj}\subset [2h_1-2,2h_1-2+b].
\]

Note that ${\rm proj}(0,1)=(2h_1-2,0)$ and that in $A_{\rm proj}$ we have a relation, namely
\begin{equation}
(0,0)+(2h_1-2,0)=2(h_1-1,0),
\end{equation}
which we did not have for the corresponding points of $A$:
\[
(0,0)+(0,1)\neq 2(h_1-1,0).
\]

Therefore, we have
\[
V(A)=h_1+b+1
\]
and
\[
V(A_1)=2h_1-1+b,
\]
so that
\[
V(A)< V(A_1).
\]

Moreover,
\[
T(A_1)=T(A)-1.
\]

In our case we have ${\rm dim}\,A=2$: we added one new  relation,  (2.7), and so according
to Theorem 3, ${\rm dim}\,A_1=1$. Thus in our trivial case Theorem 2 is proved. 
\vspace{2mm}

{\bf 2.3.} Let us give simple examples that illustrate the main idea of the proof of Theorem 2 in the case when $3\le h_2\le h_1$.

{\bf Example 1.} Take the following set $A$:

{\includegraphics[height=70mm, trim=-80mm 100mm 50mm 80mm, clip=true]{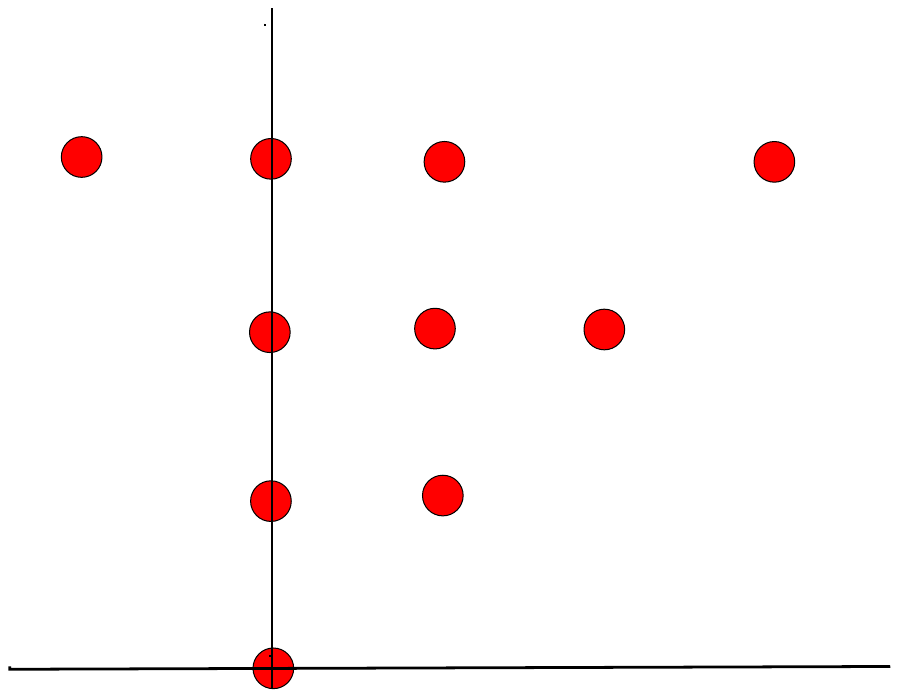}
\vspace{-30mm}
\begin{center}Figure 1\end{center}}
\[
A=\{(-1,3),(0,0),(0,1),(0,2),(0,3),(1,1),(1,2),(1,3),(2,2),(3,3)\}.
\]

Here we have $k=10$ and $T=32$.

Let us project the set $A$ onto the line $x=0$ parallel to the vector $(1,6)$. Then we have $\phi(A)=A_1$, where $\phi(-1,3)=(0,9)$, $\phi(1,1)=(0,-5)$, $\phi(1,2)=(0,-4)$, $\phi(1,3)=(0,-3)$, $\phi(2,2)=(0,-10)$, $\phi(3,3)=(0,-15)$.

\includegraphics[height=109mm, trim=-30mm 30mm 20mm 90mm, clip=true]{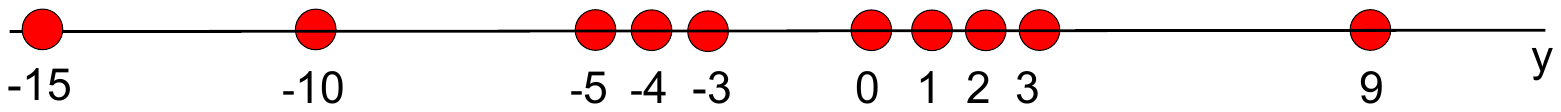}
\vspace{-90mm}
\begin{center} Figure 2\end{center}

The map $\phi$ is a one-to-one homomorphism. $\phi$ is {\it not\/} and isomorphism, because the relation $-3+3=2\cdot 0$ in $A_1$ does not yield a relation for preimages in $A$: $(0,3)+(1,3)\neq 2\cdot(0,0)$. In the present case we have
\[
T(A_1)=31,\quad V(A)=11,\quad V(A_1)=25.
\]
Therefore,
\[
V(A)<V(A_1),\quad{\rm but}\quad T(A_1)<T(A).
\]

We  have ${\rm dim}\, A=2$. Since when we pass from $A$ to $A_1=\phi(A)$ we add  one new relation (2.7), we have (see 2.1)
\[
{\rm dim}\, A_1=1.
\]

Thus, the example we just considered illustrates our theorem.

{\bf 2.4.} Let $A$ be depicted with respect to an othonormal system of coordinates. Let $\bar a=(a_1,a_2)$ and $\bar b=(b_1,b_2)$ be two points in $A$, which lie at the maximal distance (among the points of $A$) from one another (see Figure 3).

Let us draw a straight line $L_1$ through the point $\bar a$ and a straight line $L_2$ parallel to it through the point $\bar b$. We choose the lines $L_1$ and $L_2$ so that $A$ lies in the strip bounded by $L_1$ abd $L_2$, and there is a point $\bar c$ in $A$ which lies on $L_1$ or on $L_2$. Mapping the system of coordinates with the axes $L_1$ and $L_3$ into the orthonormal system of coordinates yields a map of the set $A$ onto a set $B$ which satisfies condition 1) in Section 2.4.

\includegraphics[height=120mm, trim=-95mm 0mm 0mm 0mm, clip=true]{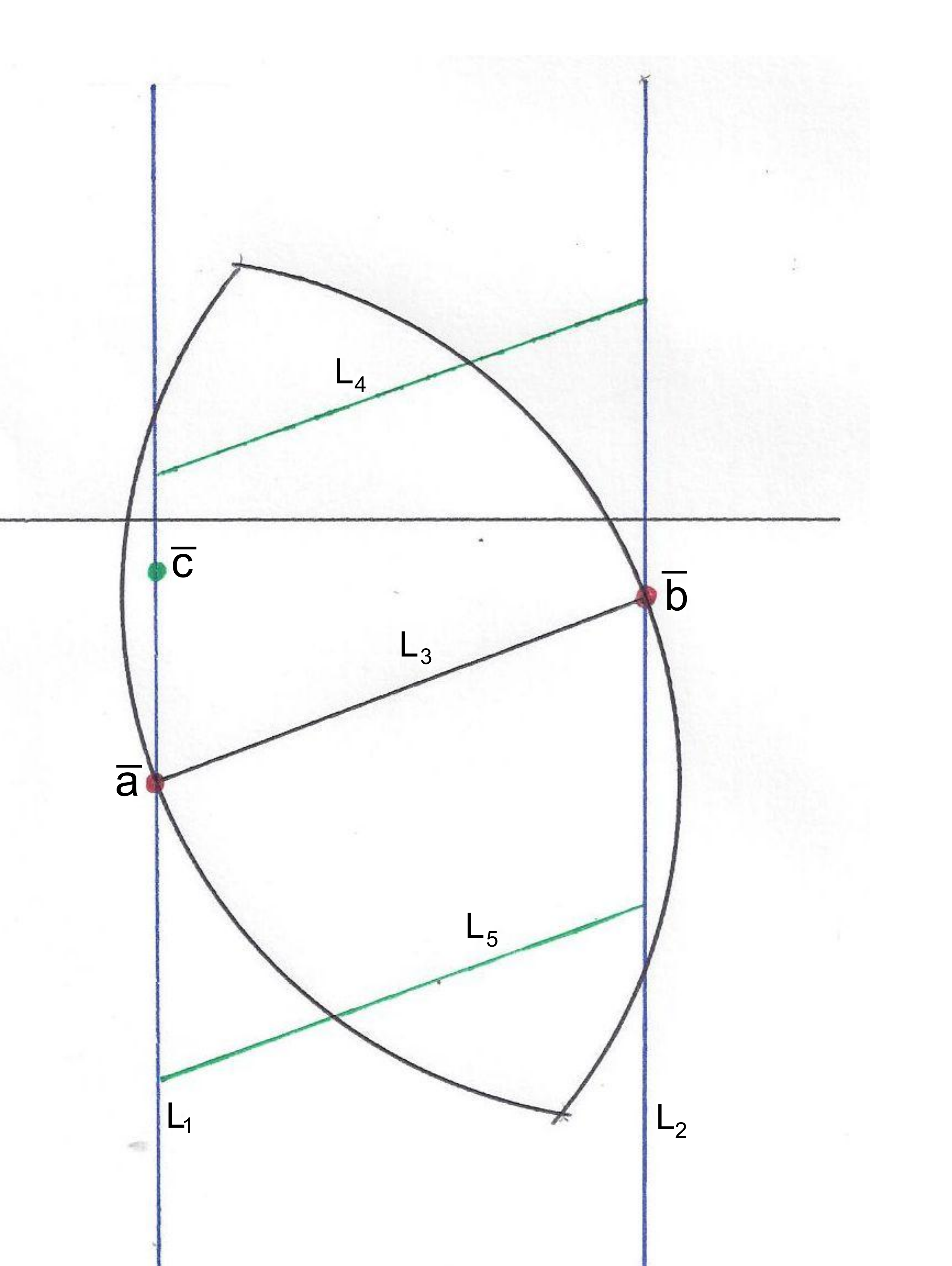}
\begin{center}Figure 3\end{center}

If we extend the example in Section 2.3 to a general case, the set $A$ lies in a strip bounded by the lines $y=0$ and
$y=h_2-1$. If $A$ is such that

 1) the points $(0,0)$ and $(0,h_2-1)$ belong to $A$, 

\noindent and

 2) the point $(1,h_2-1)$ belong to $A$,  

\noindent then passing from $A$ to $A_1=\phi(A)$ is done in the same manner as in the numerical example.

Conditions 1) and 2) in Section 2.4 cover only part of the general case presented in Section 2.7. However, here we get stronger results concerning the quantity $V (A_1)/V(A)$ and the way the presentation leads us to the core of the problem is considerably simpler that in the case treated in Section 2.7.

\includegraphics[height=100mm,trim=-50mm 0mm 0mm 76mm,clip=true]{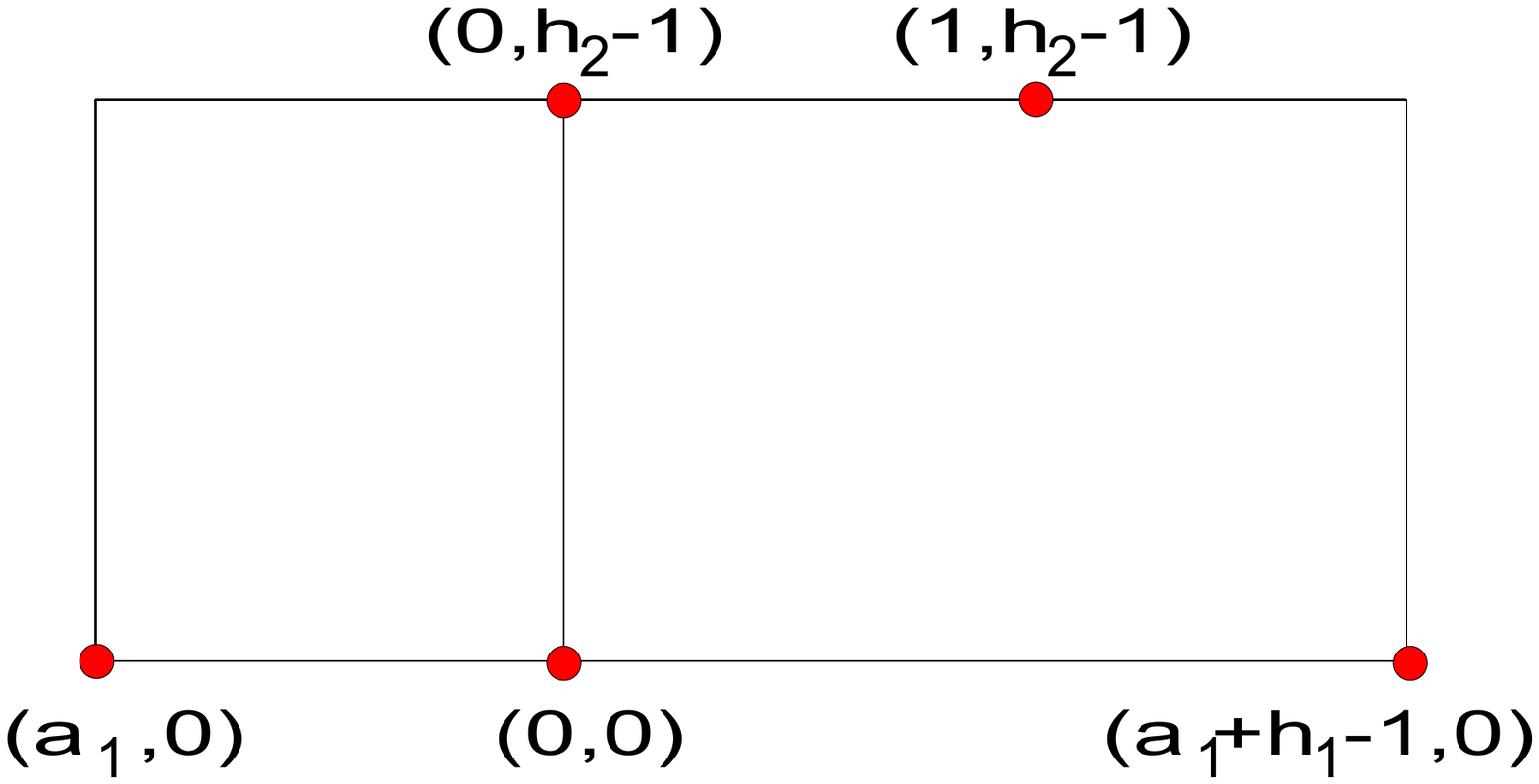}
\vspace{-50mm}
\begin{center}Figure 4\end{center}

The cases when 1)  and one of the conditions $(-1,h_2-1)\in A$, $(1,0)\in A$, or $(-1,0)\in A$ are satisfied are reduced to the case already considered by means of the reflection with respect to the line $x=0$ or to the line $y=(h_2-1)/2$. Needless to say, the shift of $A$ along the line $y=0$ does not lead to new cases.

The set $A$ with ${\rm dim}\, A=2$ is placed in the parallelepiped $D(A)$ with minimal volume $V(D)$.

According to Definition 2,
\begin{equation}
D=(a_1,a_2)+\{(x_1,x_2)\in \mathbb Z^2,\ 0\le x_1< h_1,\ 0\le x_2 <h_2\}.
\end{equation}

To satisfy the conditions 1) and 2) we put in (2.8) $a_2=0$ and take $a_1$ such that $a_1\le 0$ and $a_1+h_1-1\ge 0$.

Let us project the set $A$ onto the line $x=0$ parallel to the vector $(1,2(h_2-1))$. The matrix of this map is
 \[ 
M= \left( \begin{array}{cc}
0 & -2(h_2-1)   \\
0 & 1 \end{array} \right)
\] 
so that
\[
(x,y) \left( \begin{array}{cc}
0 & -2(h_2-1)   \\
0 & 1 \end{array} \right) =(0,y-2(h_2-1)x).
\]

On the line $x=a_1$ there are points of the set $A$. Let $(a_1,y_1)$ be  the point  with the maximal value of $y$ among those points.

In the same way we see that on the line $x=a_1+h_1-1$ there are points of $A$. Let $(a_1+h_1-1,y_0)$ be the point with the minimal value of $y$ among those points.

The images of these two points under the projection are
\[
(a_1,y_1) \left( \begin{array}{cc}
0 & -2(h_2-1)   \\
0 & 1 \end{array} \right) =(0,y_1-2(h_2-1)a_1) 
\]
and
\[
(a_1+h_1-1,y_0) \left( \begin{array}{cc}
0 & -2(h_2-1)   \\
0 & 1 \end{array} \right) =(0,y_0-2(h_2-1)(a_1+h_1-1)),
\]
respectively.

The projected set $A_1=A_{\rm proj}$ on the line $x=0$ is such that the ordinates $y$ of its points satisfy
\[
y_0-2(h_2-1)(a_1+h_1-1)\le y\le y_1-2(h_2-1)a_1.
\]

The length of this segment is equal to
\[
y_1-2(h_2-1)a_1-y_0+2(h_2-1)(a_1+h_1-1)+1=y_1-y_0+1+2(h_2-1)(h_1-1).
\]
The segment is minimal when $y_1=0$ and $y_0=h_2-1$, and then the minimal length is
\[
2(h_2-1)(h_1-1)-(h_2-1)+1=(2h_1-1)(h_2-1)+1.
\]

There is no arithmetic progression with difference $d\ge 2$ that contains the set $A_1$; indeed, otherwise $D(A)$ would not be minimal.

We conclude that
\[
V(A_1)\ge (2h_1-1)(h_2-1)+1>h_1h_2=V(A).
\]

As in Example 1, we see immediately that
\[
T(A)>T(A_1)
\]
and
\[
{\rm dim}\, A_1=1.
\]

In our case the value $V(A_1)$ is maximal compared with other possible cases.
\vspace{2mm}

{\bf 2.5.} Suppose that condition 1) of Section 2.4 is fulfilled, and 2) $(1,y_0)\in A$,
where $y_0$ is maximal on the line $x=1$ and
\begin{equation*}
y_0\ge \frac{1}{2}(h_2-1).\tag{$\ast$}
\end{equation*}

The case when $y_0\le \frac{1}{2}(h_2-1)$ can be studied by using reflection with respect to the line $y= \frac{1}{2}(h_2-1)$. The case when $y_0=h_2-1$, was studied in Subsection 2.4, is of course a particular case of the case studied below in Subsection 2.5, but we provided it in order to illustrate the proof of Theorem 2 on a simpel example, which, incidentally, also gives the bigger value of $V(A_1)$.  

Let us first illustrate our general case by means of the following example.
\vspace{2mm}

{\bf Example 2.} Consider the set (see Figure 5).

\includegraphics[height=160mm, trim=-30mm  0mm 0mm 90mm, clip=true]{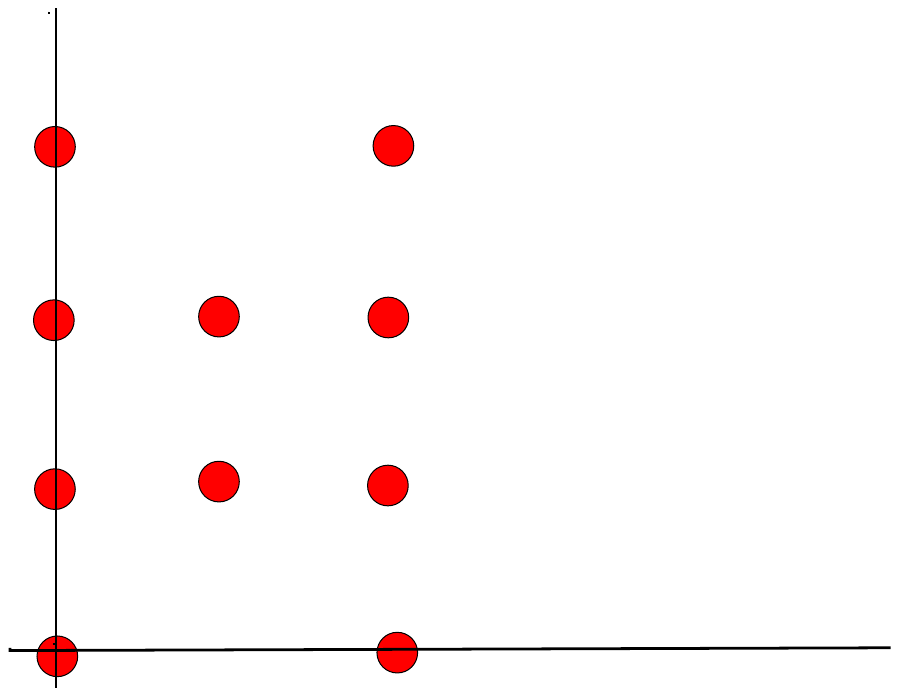}
\vspace{-110mm}
\begin{center}{Figure 5}\end{center}

\[
A=\{(0,0),(0,1),(0,2), (0,3),(1,1),(1,2),(2,0),(2,1),(2,2),(2,3)\}.
\]

In this case $k=10$ and $T=33$.

Let us project $A$ onto the line $x=0$ parallel to the vector $(1,5)$. We obtain the set $A_1=\phi(A)$, where $\phi(1,2)=(0,-3)$, $\phi(1,1)=(0,-4)$, $\phi(2,3)=(0,-7)$, $\phi(2,2)=(0,-8)$, $\phi(2,1)=(0,-9)$, $\phi(2,0)=(0,-10)$. The points on the line $x=0$ remain on place. 

\includegraphics[height=95mm, trim=80mm 50mm 0mm 160mm, clip=true]{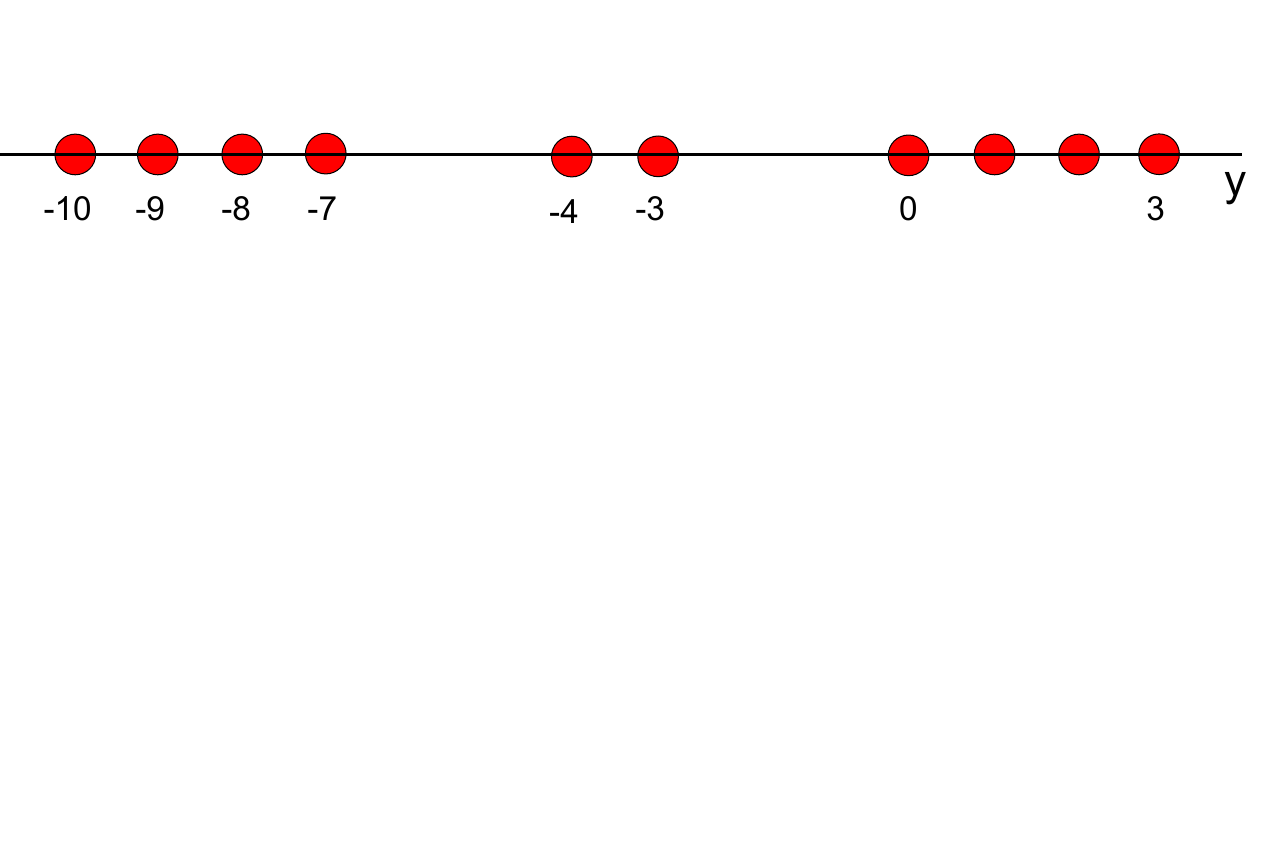}
\vspace{-60mm}
\begin{center}{Figure 6}\end{center}

We have $V(A)=12$ and $V(A_1)=14$. Moreover, $T(A)=33$, while $T(A_1)=32$. Therefore,
\[
V(A)<V(A_1)
\]
and
\[
T(A)>T(A_1).
\]

\newpage

{\bf 2.6.} Now let $\ell=(1,y_1)$, where $y_1=h_2-1+y_0$.

Let us project the set $A$ onto the line $x=0$ parallel to the vector $\ell$. Let $y_0$ be the maximal number characterizing the points in conditions 1) and 2) of Subsection 2.5 (this condition may be satisfied by using a suitable isomorophic transformation; it give the set $A_1$ with maximal volume $V(A_1)$.

The projection in question is defined by the matrix
\[
M=\left( \begin{array}{cc}
0 & -(h_2-1+y_0)   \\
0 & 1 \end{array} \right), 
\]
so that
\[
\phi(x,y)=(x,y)\left( \begin{array}{cc}
0 & -(h_2-1+y_0)   \\
0 & 1 \end{array} \right) = (0,y-(h_2-1+y_0)x).
\]

Let
\[
I_a=\{(x,y)|\, x=a,\ 0\le y\le h_2-1\}.
\]
Then
\[
\phi(I_a)=\{(x,y)|\, x=0,\ -(h_2-1+y_0)a\le y\le h_2-1-(h_2-1+y_0)a\}.
\]
Thus we can write
\[
\phi(D)=\bigcup_a \phi(I_a).
\]

Now let us consider, in the same way as in Subsection 2.4, the points $(a_1,y_1)$ and $(a_1+h_1-1,y_2)$. The images of these two points under $\phi$ are
\[
(a_1,y_1)M=(0,y-a_1(h_2-1+y_0))
\]
and
\[
(a_1+h_1-1,y_2)M=(0,y_2-(a_1+h_1-1)(h_2-1+y_0)),
\]
respectively.

The projected set $A_1$ on the line $x=0$ is such that for $y\in A_1$ we have
\[
y_1-a_1(h_2-1+y_0)\le y\le y_2-(a_1+h_1-1)(h-2-1+y_0).
\]
The length of this segment is
\[
y_2-(a_1+h_1-1)(h_2-1+y_0)-y_1+a_1(h_2-1+y_0)+1=y_2-y_1+(h_2-1+y_0)(h_1-1).
\]
This segment has a minimal length when $y_1=0$ and $y_2=h_2-1$, and then its length is
\[
(h_2-1+y_0)(h_1-1)-(h_2-1)+1.
\]

We conclude that it holds that
\begin{equation*}
V(A_1)\ge (h_2-1+y_0)(h_1-1)-(h_2-1)+1>h_1h_2=V(A).\tag{$\ast\ast$}
\end{equation*}

In view of ($\ast$), we have
\[
(h_2-1+y_0)(h_1-1)-(h_2-1)+1\ge (h_2-1)(h_1-2)+\frac{h_2-1}{2}(h_1-1)+1=
\]
\[
=h_1h_2-h_1-2h_2+3+\frac{h_2-1}{2}(h_1-1).
\]
To get ($\ast\ast$) we simply verify that
\[
\frac{(h_2-1)(h_1-1)}{2}-h_1-2h_2-3>0,
\]
or
\[
\frac{(h_2-3)(h_1-1)}{2}-2h_2+2>0,
\]
or, further
\[
(h_2-3)(h_1-1)-4(h_2-3)-8>0
\]
or
\[
(h-3)(h_1-5)>8.
\]

If $h_2>3$ and $h_1>13$ this last inequality is satisfied, and the last conditions are fulfilled when $k$  is sufficiently large and $h_2>3$. 
\vspace{2mm}
\indent The analysis of the case $h_2=3$ is left as an exercise.
\vspace{2mm}

{\bf 2.7.} Suppose now that condition 1) of  Section 2.4 is fulfilled and
\[
{\rm 2)}\qquad (x_0,y_0)\in A,
\]
where $|x_0|$ is  minimal among the points $(x,y)\in A$, for which $x\neq 0$ and
\[
|x_0|>1.
\]
We can assume (using reflexion) that
\[
x_0>0.
\]

Let us take all the $\{x\}$-coordinates of points of $A$ such that
\[
a_1\le x\le a_1+h_1-1
\]
and $a_1\le 0$, and label them as
\[
x_1=a_1,\,x_2,\dots x_r=a_1+h_1-1.
\]

For the points $(x_i,y)\in A$, denote $\max_{(x_i,y)\in A}y=y_{i\,{\rm max}}$ and $\min_{(x_i,y)\in A}y=y_{i\,{\rm min}}$.

Consider also the points
\[
(x_{i+1},y_{i+1\,{\rm min}})\quad{\rm and}\quad (x_{i+1},y_{i+1\,{\rm max}})
\]

\includegraphics[height=120mm, trim=-50mm 0mm 0mm 45mm, clip=true]{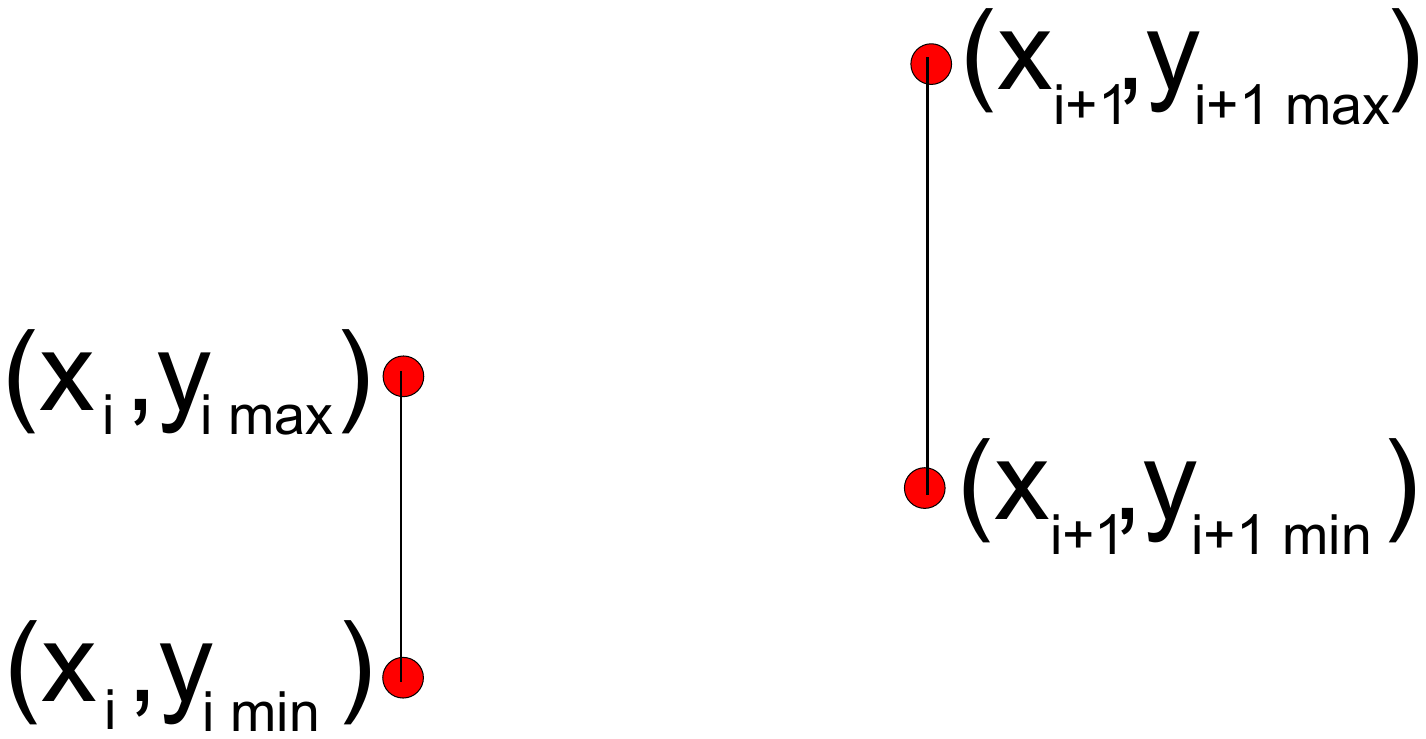}

\vspace{-60mm}

Further, let us denote
\begin{equation}
\Delta_i^x:= x_{i+1}-x_i 
\end{equation}
and
\begin{equation}
\Delta_i^y:=\max\left\{y_{i\,{\rm max}}-y_{i+1\,{\rm min}}, y_{i+1\,{\rm max}}-y_{i\,{\rm min}}\right\}.
\end{equation}

Note that by using, if necessary, the reflection of $A$ with respect to the line $y=(h_2-1)/2$, we in fact have that
\begin{equation}
\Delta_i^y:= y_{i+1\,{\rm max}}-y_{i\,{\rm min}}, 
\end{equation}
and we will assume this from now on.

Further, let us denote
\begin{equation}
\tau = \min_i\frac{\Delta_i^x}{(\Delta_i^y+h_2-1)} 
\end{equation}

In this section we will consider only the case when there exists an $i$ such that
\begin{equation}
\Delta_i^x=1.
\end{equation}

Needless to say, there may be several values of $i$ such that (2.13) holds. We denote by $i_0$ that value of $i$  which gives the minimum in $\tau$.

Now let us consider the projection $\phi$ along the vector
\[
\ell:=(1,\Delta_{i_0}^y+h_2-1).
\]
onto the line $x=0$ (in the definition of $\ell$ we take $y_1=\Delta_{i_0}^y+h_2-1$ instead of $y_1=h_2-1+y_0$, as in Subsection 2.6). 

The proof of the inequality $V(A_1)\ge V(A)$ is completed as in Subsection 2.6, under the additional condition 
\begin{equation}
\Delta_{i_0}^y\ge 1.
\end{equation}

Note that we have
\[
\phi(x+1,y_{x+1\,{\rm max}})=\phi(x,y_{x+1\,{\rm max}}-(h_2-1)-(y_{x+1\,{\rm max}}-y_{x\,{\rm max}}))=\phi(x,y_{x\,{\rm max}}-(h_2-1)).
\]
The points $(x,y_{x\,{\rm max}})$ and $(x,y_{x\,{\rm max}}-(h_2-1))$ give for ordinates the same difference as for the original points $(0,0)$ and $(0,h_2-1)$, i.e., in $A_1$ we have a new relation compared with $A$, and
\[
{\rm dim}\,A_1=1.
\]
\vspace{2mm}

{\bf 2.8.}

Let us explain this connection in more detail.

Define, using the notations in Section 2.6,
\[
H_a=\left\{(x,y)|\, x=0,\,\, h_2-1-(h_2-1+y_0)a\le y\le (h_2-1+y_0)(a+1)\right\},
\]
for 
\[
a_1\le a\le a_1+h_1-2.
\]

The points in the set $A_1''=\bigcup_{a=a_1}^{a_1+h_1-2} H_a$ have no preimages under the map $\phi$, and the points of the set $A_1'=A_1\setminus  A_1''$ have such preimages, and this fact explains why $V(A)\le V(A_1)$.

In the case considered in our Subsection 2.7  we have a corresponding set $B_1''$ which has no preimages, and which is defined in the following manner:
\[
B_1''=\bigcup_{a=a_1}^{a_1+h_1-2} F_a,
\]
where
\begin{equation}
F_a=\left\{(x,y)|\, x=0,\,\, h_2-1-(h_2-1+\Delta_{i_0}^y)a\le y\le (h_2-1+\Delta_{i_0}^y)(a+1)\right\}.
\end{equation}

In order for the set of admissible $y$ in (2.15) to be non-empty it is necessary that either (2.14) holds or $\Delta^y_{i_0}=0$. In the case where (2.14) holds we have $V(A_1)\ge V(A)$.  
\vspace{2mm}

{\bf 2.9.} Let us give an example for the case considered in Subsection 2.8.

{\bf Example.} Consider the set
\[
A=\{(0,0),(0,1),(0,2),(0,3),(0,4),(2,2),(3,2),(3,3),(4,2),(5,2)\}
\]
or, graphically
\newpage
\includegraphics[height=120mm,trim=-50mm 30mm 50mm 70mm,clip=true]{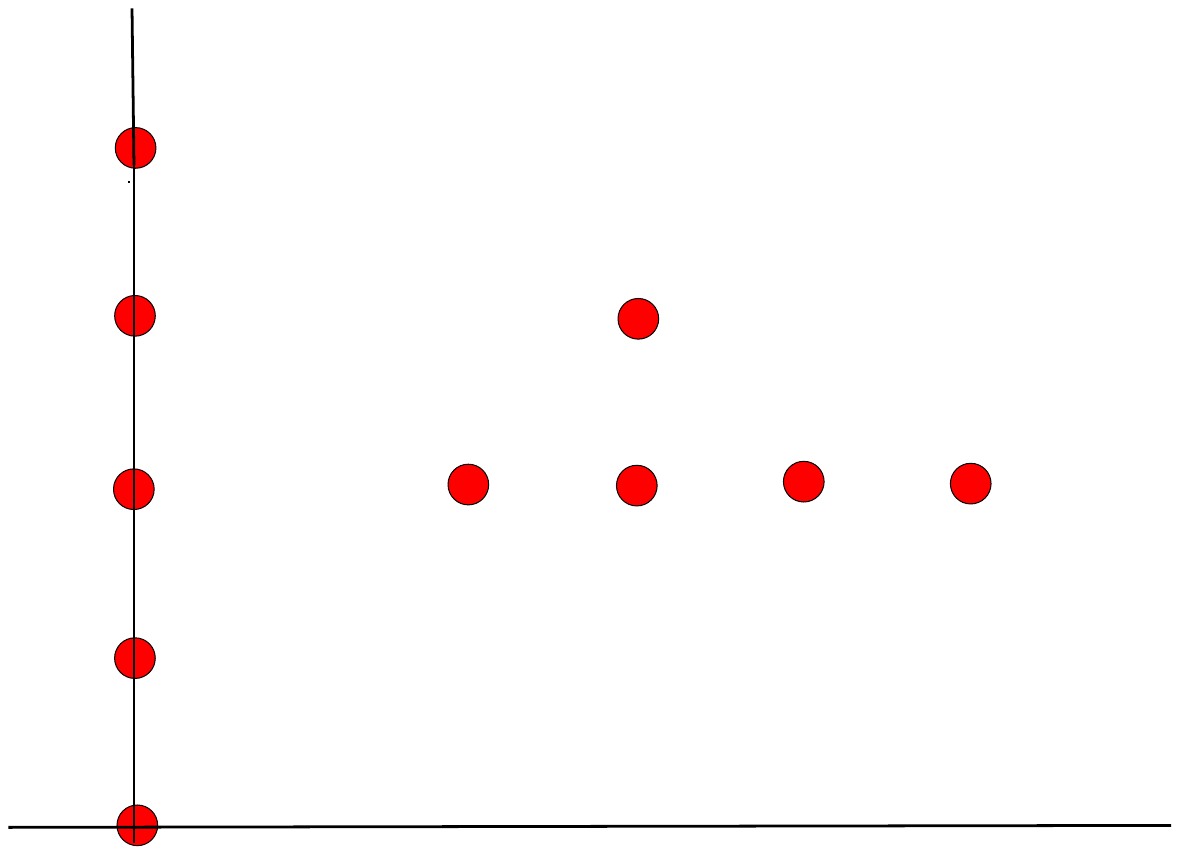}
\vspace{-50mm}

Taking the projection $\phi$ onto the line $y=0$ along the vector $(1,5)$ we obtain the set
\[
A_1=\phi(A)=\{(0,0),(0,1),(0,2),(0,3),(0,4),(0,-8),(0,-12),(0,-13),(0,-18),(0,-23)\}.
\]

We have
\[
V(A)\le 5\cdot 6-2=28 =V(A_1)
\]
and
\[
T(A)=25,\qquad T(A_1)=24,\qquad{\rm i.e.,}\,\, T(A_1)<T(A).
\]
Thus, Theorem 2 is proved. 

\vspace{2mm}

{\bf 2.10.} When we are studying the case where $\Delta^y_{i_0}=0$, we chose the  vector $\ell$ defining the projection without asking that the value $|x_0|$ be minimal; rather, we just keep the requirement that the projection be one-to-one.

We explain this with the help of the following example. Consider the set
\[
A=\{(0,0),(0,1),(0,2),(0,4),(0,8),(3,8),(4,8),(5,8),(6,8),(9,8),(9,4),(9,2),(9,1),(9,0)\}
\]

Here $V(A)=90$.

Let us contract the set $A$ along the $x$-axis: $\phi_1(x,y)=(x/3,y)$:
\[
A_1=\phi_1(A)=\{(0,0),(0,1),(0,2),(0,4),(0,8),(1,8),
\]
\[
({4}/{3},8),({5}/{3},8),(2,8),(3,8),(3,4),(3,2),(3,1),(3,0)\}
\]

Now let us project the set $A$ onto the line by means of the vector $(1,16)$. We obtain:
\[
A_2=\phi_2(A_1)=\{(0,8),(0,4),(0,2),(0,1),(0,0),(0,-8),
\]
\[
(0,-13\cdot {1}/{3}), (0,-18\cdot {2}/{3})(0,-24), (0,-40),(0,-44),(0,-46)(0,-47), (0,-48)\}. 
\]

Finally, let us dilate the set $A_2$ along the $y$ axis:   $\phi_3(x,y)=(x,3y)$. We get
$
A_3=\phi_3(A_3)$
\[
=\{(x,y)|\, x=0,\,  y=24,12,6,3,0,-24,-40,
-56,-72,-120,-132,-138,-141,-144\}.
\]

Therefore,
\[
V(A_3)=169
\]
and 
\[
V(A_3)>V(A).
\]

In $A_3$ we have
\[
24+(-24)=2\cdot 0=0,
\]
but for the preimages in $A$ we have
\[
(0,8)+(3,8)\neq =0,
\]
so that ${\rm dim}\, A_3=1$.

The transition from dimension $d+1$, with $d\ge 2$, to dimension $d$, is carried out in similar manner.

\vspace{2mm}

\section{On $L_{\rm m}$, the exact value of the extremal volume for $\mathbb A(k,T)$ }

{\bf 3.1.  Length of extremal sets.} 

In Section 3 we will study sets for which ${\rm dim}\, A=1$. Such sets lie inside segments of $\mathbb Z$. The minimal length of a family of
isomorphic sets is given by a set in normal form which lies in the segment $[0, a_{k-1}]$. Let us define ${\mathbb A}(k,T)$ as the set of all families of
isomorphic sets with given values of $k$ and $T$. Two sets in ${\mathbb A}(k,T)$ are not necessarily isomorphic. Denote
\[
a_{\rm m}=\max_{A\in {\mathbb A}(k,T)} a_{k-1}.
\]

A set for which $a_{k-1}=a_{\rm m}$ will be called {\it extremal\/}. The volume (i.e., length) of such a set is equal to
\[
V(A)=L(A)=1+a_{\rm m}.
\]

\indent{\bf Exact formula for $a_{\rm m}$.}

In [4], p. 4 the following conjecture was formulated:

{\it If
\begin{equation}
T=ck-\frac{c^2+c-4}{2}+b,
\end{equation}
where
\begin{equation}
2\le c\le k-1
\end{equation}
and
\begin{equation}
0\le b\le k-c-1,
\end{equation}
then
\begin{equation}
a_{\rm m}=2^{c-2}(k+1-c+b).
\end{equation}}

{\bf 3.2. The $T\longleftrightarrow (c,b)$ one-to-one correspondence}
 
We will ensure that the function $a_{\rm m}={\cal F}(k,T)$ is well defined if we can show that the set of values $\{T\}$ and $\{(c,b)\}$ are into a one-to-one
correspondence.

For a given $c$, the values of $T$ are, in view of (3.1) and (3.2), situated in the segment
\[
\left[ck-\frac{c^2+c-4}{2},ck-\frac{c^2+c-4}{2}+k-c-1\right]=
\]
\[
=\left[ck-\frac{c^2+c-4}{2},(c+1)k-\frac{(c+1)^2+(c+1)-4}{2}\right].
\]
$T$ is minimal when $c=2$ and $b=0$, i.e.,
\[
T=2k-1.
\]

$T$ is maximal when $c=k-1$ and $b=0$, i.e.,
\[
T=(k-1)k- \frac{(k-1)^2+(k-1)-4}{2}= \frac{k^2-k}{2}+2.
\]

{\bf 3.3. Concrete example.} Let us give a concrete numerical example explaining how to obtain an estimate from below for ${\cal F}(k,T)$ (we will conjecture that this estimate is
in fact sharp).

Let us take the set
\[
A=\{0,1,2,4,8,16,32,48,96,192,384\}.
\]
Here we have $|A|=11$ and $a_{\rm m}=384$. We will calculate $T$ in several steps, appropriate for further generalization.

Denote
\[
A_1=\{8,16,32,48\}=8\times\{0,1,2,4,6\},
\]
\[
A_2=\{1,2,4\}=\{2^0,2^1,2^2\},
\]
and
\[
A_3=\{96,192,384\}=96\times \{2^0,2^1,2^2\}.
\]

The sums $1+a_i$ contribute 10 numbers to $T$, $A\setminus\{1\}$ is left and $2+a_i$ contributes 9 sums, the number 4  gives 8 sums, 384 gives 7 sums, 192
gives 6 sums, and 96 gives 5 sums. We get the set  $8\times\{0,1,2,4,6\}$. Here $k'= 5$, $b=2$ (the number of holes), and $T'=10-1+2=11$. Consequently,
\[
T=10+9+8+7+6+5+11=56.
\]

{\bf 3.4. General example.} Now it is not difficult, using the above example, to build a generalized example $A$. To this end we will use  numbers $k_0=k_2+k_3$, $k_1$ and $b\le k_1-3$.
We take the set $A$ in the form
\[
A=A_1\cup A_2\cup A_3,
\]
where
\[
A_2=\{1,2,2^2,\ldots,2^{k_2-1}\},\qquad A_3=2^{k_2}\times A_3'.
\]

The set $A_1'$ is defined by the two numbers $k_1$ and $b$ by means of the Basic Theorem.  Namely, $A'_1$ is a part of the segment $[0,k_1-1+b]$ satisfying
the condition $|2A_1'|=2k_1-1+b$. The detailed description of the set $A'_1$ is provided by the Basic Theorem.

The last number in the set $A_1$ is $p=2^{k_2}(k_1-1+b)$.                               
Now let us describe the set $A_3$:
\[
A_3=2p\times\{1,2,2^2,\ldots,2^{k_3-1}\}.
\]
Put $k_0=k_2+k_3$. Computing $T$ as in the specific example above, we get
\[
T=k-1+(k-2)+\cdots+(k-k_2)+(k-k_2-1)+\cdots+(k-k_0)+2k_1-1+b=
\]
\[
= \frac{k-k_0+k-1}{2}\,k_0+2k_1-1+b=
\]
\[
=
k_0k-\frac{(k_0+1)k_0}{2}+2k-2k_0+b-1= (k_0+2)k-\frac{(k_0+5)k_0}{2}+b-1.
\]
Denoting $k_0+2+c$, we get
\[
T=ck-\frac{c^2+c-4}{2}+b,
\]
as in (3.1).

Here
\[
2\le c\le k-1,\qquad 0\le b\le k-c-1.
\]

{\bf 3.5.}   Not very much is known about the exact form of the function ${\cal F}(T)=a_{\rm m}$ (only, the conjecture below).

The Basic Theorem gives us:

If $2k-1\le T\le 3k-4$, when
\[
T=2k-1+b,\quad 0\le b\le k-3,
\]
then
\[
a_{\rm m}=k-1+b.
\]

Strangely enough, even for
\[
T=3k-4+b,\quad 0\le b\le k-4,
\]
in which case, according to my early conjecture, it must hold that
\[
a_{\rm m}=2k-4+2b 
\]
the proof is still missing.

In a series of papers, Renling Jin (see [8], [9], [10])  studied  similar problems.

\section{Approximate groups}  

B. Green and T. Tao have introduced, in a series of seminal papers, the notion of approximate group.
For given $k$ and $T$, denote by ${\mathbb {AG}}(k,T)$ the collection of all sets that are {\it approximate groups in\/} 
$\mathbb Z$ (Tao, Green) and have characteristic values  $k$ and $T$.

In our case such sets have $A$ have a simple description: $A$ is a finite set of integers with characteristics $k$ and $T$ (as everywhere in the paper) and satisfy the additionaal conditions

(1) $A$ is symmetric, i.e., if $a\in A$, then $-a\in A$;

(2) $0\in A$.

The set ${\mathbb{ AG}}(k,T)$ is of course partitioned into classes with respect to the relation given by isomorphism.

As it was explained in Section 2, we can focus on the study of sets $A$ with ${\rm dim}\, A=1$. In this case the volume (or legth) of a set $A$ can be calculates as follows:

We may assume that ${\rm g.c.d.}(A)=1$. The there exists a segment 
$L_{\rm m}=[-b_{\rm m},b_{\rm m}]$ such that $A\subset L_{\rm m}$ and $L_{\rm m}$ has minimal length. Then
\[
V(A)=V(L_{\rm m})=2b_m+1.
\]

We propose the following 

{\bf Conjecture.} {\it For a  set $A\in {\mathbb {AG}}(k,T)$ with $T$ written in the form
\[
T=ck-\frac{\displaystyle{3c^2-2c-4}}{\displaystyle{2}}+b
\]
we have
\begin{equation}
L_{\rm m}=3^{\frac{c}{2}-1}(k-c+b+1)+1,
\end{equation}
with
\[
2\le c\le k-1,\qquad 0\le b\le k-c-1.
\]}

Let us construct a set for which $L_{\rm }$ is given by the right-hand side of (3.5 ). This example will show that $L_{\rm m}$ cannot be smaller, i.e.,
\[
L_{\rm m}\ge 3^{\frac{c}{2}-1}(k-c+b+1)+1.
\]

Let $|A|=k$, $k$ odd, and assume $A$ is symmetric with respect to 0. Suppose the first $\overline k_1$ elements, $|A_0|=\overline k_1$, obey the Basic Theorem. This means that ($b$ is even)
\begin{equation}
T_0=2\overline k_1-1+b\qquad L=\overline k_1+b.
\end{equation}

$A_0$ lies in the following segment:
\[
A\subset\left[-\frac{\overline k_1-1+b_1}{2},\frac{\overline k_1-1+b_1}{2}\right]=[p,p].
\]

Let us add here two points, $-3p$ and $3p$, then another two points $-9p$ and $9p$, and so on. Overall, we add $2\overline k_2$ points, so we have the segment
\[
\left[-3^{\overline k_2}p,3^{\overline k_2}p\right].
\]
 
We have a total of $k$ points, with
\[
k=\overline k_1+2\overline k_2.
\]

The set $A$ lies in a segmeny of length  $L$ (i.e., the segment contains $L$ points), where
\[
L=3^{\overline k_2}\cdot 2p+1 = 3^{\overline k_2}\cdot 2\,\frac{\overline k_1-1+b_1}{2}+1=
\]
\[
=3^{\overline k_2}\left(\overline k_1+b-1\right)+1=3^{\overline k_2}\left(k-2\overline k_2+b-1\right)+1.
\]
Denoting $c = 2\overline k_2+2$, we have
\[
L=3^{\frac{c}{2}-1}(k-c+b+1)+1.
\] 

From the Basic Theorem it follows that
\[
0\le b\le \overline k_1-3.
\]
 
\qquad Now let us calculate $T$. 
For the first $\overline k_1$ points we assume that the Basic Theorem holds true, so that
\[
T_0=2\overline k_1-1+b, \qquad  0\le b\le \overline k_1-3
\]
($b$ is the number of holes, $b/2$ holes to the left, and $b/2$ to the right). 
Let us add two more points (see (3.6)); then to $T_0$ one adds $2\overline k_1$. Indeed, after we add the pair  of points 
\[
-3p,\dots,-p,0,\dots,p,3p,
\]
we have $3p+(-p)=2p$, $3p+(-3p)=0$. From $-p$ to $p$ there are $\overline k_1$ points: $3p$ with these points, and $-3p$ of these points contribute $2\overline k_1$ to $T_0$. Adding two more points yields
\[
2\overline k_1+2=2\left(\overline k_1+1\right).
\]  
Further, we get
\[
2\overline k_1=2\left(\overline k_1+2\right).
\]
 
Let us a pair $\overline k_2$ times; the additions to $T$ are given by
\[
2\overline k_1,\, 2\left(\overline k_1+1\right),\,2\left(\overline k_1+2\right),\,\dots\,2\left(\overline  k_1+\overline k_2-1\right).
\]
The last pair was $k$.

Altogether we get
\[
T=T_0+2\overline k_1+2\left(\overline k_1+1\right)+\cdots+2\left(\overline  k_1+\overline k_2-1\right)=
\]
\[
=T_0+2\overline k_2\frac{\overline k_1+\overline k_1+\overline k_2-1}{2}= 
\]
Here there are $\overline k_2$ terms (in an arithmetic progression). The first term is $2\overline k_1$, and the last is $2(\overline k_1+\overline k_2-1)$, so that continuing we have
\[
=T_0+\overline k_2(2\overline k_1+\overline k_2-1)= T_0+\overline k_2(2k-3\overline k_2-1)=2\overline k_1=1+b+\cdots \, .
\]

We have
\[
\overline k_2=\frac{2\overline k_2+2}{2}-1= \frac{c}{2}-1,
\]
so that
\[
T=2k-4\overline k_2-1+b+\left(\frac{c}{2}-1\right)\left(2k-\frac{3c}{2} +3-1\right)=
\]
\[
=2k-2c+4-1+b+ck-2k-\frac{3c^2}{4}+\frac{3c}{2}+c-2=
\]
\[
=ck+\frac{c}{2}-\frac{3c^2}{4}+1+b=ck-\frac{3c^2-2c-4}{4}+b.
\]
This is in agreement with the above conjecture.

\section{} In conclusion, let us note that the description of the exact structure of extremal sets in the Abelian case will clearly lead to similar progress in the non-Abelian case. The papers [9], [10], and [11] give some hints on how to begin such a study.
\vspace{5mm}

{\bf References}

[1]\quad  Y.  Bilu, ``\emph{Structure of sets with small sumset}'',
Ast\'erisque {\bf 258} (1999), 77--108.
\vspace{2mm}

[2]\quad  Chang, Mei Chu, ``\emph{A polynomial bound in Freiman's theorem}'',
Duke Math. J. {\bf 113} (2002), no. 3, 399--419.
\vspace{2mm}

[3]\quad G. A. Freiman, ``\emph{What is the structure of $K$ if $K+K$ is small?}'', Lecture Notes in Math. 1240, 109--134, Springer, 1987.
\vspace{2mm}

[4]\quad Freiman, Gregory A. ``\emph{Structure theory of set addition {\rm  III}. Results and problems}". (2012)  arXiv:1204.5288.
\vspace{2mm}

[5]\quad  Freiman, Gregory, Herzog, Marcel, Longobardi, Patrizia, Maj, Mercede, ``\emph{Small doubling in ordered groups}''. J. Aust. Math. Soc. 96 (2014), no. 3, 316--325.
\vspace{2mm}

[6]\quad Freiman, G. A., Herzog, M., Longobardi, P., Maj, M., Stanchescu, Y. V., ``\emph{Direct and inverse problems in additive number theory and in non-abelian group theory}''. European J. Combin. 40 (2014), 42--54.
\vspace{2mm}

[7]\quad  Freiman, G. A., Herzog, M., Longobardi, P., Maj, M., Stanchescu, Y. V., ``\emph{A small doubling structure theorem in a Baumslag–Solitar group}''. European J. Combin. 44 (2015), part A, 106--124.
\vspace{2mm}

[8]\quad Jin, Renling, ``\emph{Freiman's inverse problem with small doubling property}''.
Adv. Math. 216 (2007), no. 2, 711-752.
\vspace{2mm}

[9]\quad Jin, Renling, ``\emph{Characterizing the structure of A when the ratio $|2A|/|A|$ is bounded by 3+epsilon}''.   arXiv:math/0504186. 
\vspace{2mm}

[10]\quad  Renling Jin, ``\emph{Detailed structure for Freiman's $3k-3$ Theorem}''. arXiv:1308.0741.
\vspace{2mm}

[11]\quad  S. V. Konyagin, ``\emph{On Freiman's Theorem}'', in: Abstracts of Paul Turan Memorial Conference, Budapest, Hungary, August 22--26, 2011.
\vspace{2mm}

[12]\quad Konyagin, Sergei V., Lev, Vsevolod F., ``\emph{Combinatorics and linear algebra of Freiman's isomorphism}''. Mathematika 47 (2000), no. 1-2, 39--51
\vspace{2mm}

[13]\quad M. Nathanson, ``\emph{Additive Number Theory. Inverse Problems and the Geometry of Sumsets}'', Graduate Texts in Math., Vol. 165, Springer Verlag, 1996.
\vspace{2mm}

[14]\quad  I. Ruzsa, ``\emph{Generalized arithmetical progressions and sumsets}'',
Acta Math. Hungar. {\bf 65} (1994), no. 4, 379--388.
\vspace{2mm}

[15]\quad T. Sanders, ``\emph{Appendix to: `Roth's theorem on progressions revisited}'',  \text{[J. Anal. Math.} \text{ {\bf 104} (2008), 155--192] by J. Bourgain},
J. Anal. Math. {\bf 104} (2008), 193--206.
\vspace{2mm}

[16]\quad T. Schoen, ``\emph{Near optimal bounds in Freiman's theorem}'', Duke Math. J. {\bf 158} (2011), no. 1,  1--12.


\end{document}